\input amstex
\input amsppt.sty
\magnification=\magstep1
\vsize=22.2truecm
\baselineskip=16truept
\nologo
\pageno=1
\topmatter
\def\Z{\Bbb Z}
\def\N{\Bbb N}
\def\K{\Bbb K}
\def\Q{\Bbb Q}
\def\CO{\Cal O}
\def\fp{\frak p}
\def\fP{\frak P}
\def\Proof{\noindent{\it Proof}}

\def\Ack{\medskip\noindent {\bf Acknowledgment}}
\def\pmod #1{\ (\roman{mod}\ #1)}

\def\jacob #1#2{\left(\frac{#1}{#2}\right)}
\title
A generalization of Wolstenholme's harmonic
series congruence
\endtitle
\author
Hao Pan
\endauthor
\abstract Let $A, B$ be two non-zero integers. Define the Lucas sequences $\{u_n\}_{n=0}^{\infty}$ and $\{v_n\}_{n=0}^{\infty}$ by
$$
u_0=0,\ u_1=1,\ u_n=Au_{n-1}-Bu_{n-2} \text{ for } n\geq 2
$$
and
$$
v_0=2,\ v_1=A,\ v_n=Av_{n-1}-Bv_{n-2} \text{ for } n\geq 2.
$$
For any $n\in\Z^+$, let $w_n$ be the largest divisor of $u_n$ prime to $u_1,u_2,\ldots,u_{n-1}$.
We prove that for any $n\geq 5$
$$
\sum_{j=1}^{n-1}\frac{v_j}{u_j}\equiv\frac{(n^2-1)\Delta}{6}\cdot\frac{u_n}{v_n}\pmod{w_n^2},
$$
where $\Delta=A^2-4B$.
\endabstract
\address
Department of Mathematics, Nanjing University,
Nanjing 210093, People's Republic of China
\endaddress
\email{haopan79\@yahoo.com.cn}\endemail
\endtopmatter
\document
\TagsOnRight
\heading
1. Introduction
\endheading
Let $A, B$ be two non-zero integers. Define the Lucas sequence $\{u_n\}_{n=0}^{\infty}$ by
$$
u_0=0,\ u_1=1\text{ and }u_n=Au_{n-1}-Bu_{n-2}\text{ for }n\geq 2.
$$
Also its companion sequence $\{v_n\}_{n=0}^{\infty}$ is given by
$$
v_0=2,\ v_1=A\text{ and }v_n=Av_{n-1}-Bv_{n-2}\text{ for }n\geq 2.
$$
Let $\Delta=A^2-4B$ be the discriminant of $\{u_n\}_{n=0}^{\infty}$ and $\{v_n\}_{n=0}^{\infty}$.
It is easy to show that
$$
v_n=\alpha^n+\beta^n
$$
and
$$
u_n=\sum_{j=0}^{n-1}\alpha^j\beta^{n-1-j}=\cases
n\alpha^{n-1}\qquad&\text{if
}\Delta=0,\\(\alpha^n-\beta^n)/(\alpha-\beta)\qquad&\text{otherwise},\endcases
$$
where
$$
\alpha=\frac{1}{2}(A+\sqrt{\Delta}),\ \beta=\frac{1}{2}(A-\sqrt{\Delta}).
$$

Let $p\geq 5$ be a prime. The well-known Wolstenholme's harmonic
series congruence asserts that
$$
\sum_{j=1}^{p-1}\frac{1}{j}\equiv 0\pmod{p^2}.\tag 1.1
$$
In [3], Kimball and Webb proved a
generalization of (1.1) involving the Lucas sequences. Let $r$ be the
rank of apparition of $p$ in the sequence
$\{u_n\}_{n=0}^{\infty}$, i.e., $r$ be the least positive
integer such that $p\mid u_r$. Kimball and Webb showed that
$$
\sum_{j=1}^{r-1}\frac{v_j}{u_j}\equiv 0\pmod{p^2}\tag 1.2
$$
provided that $\Delta=0$ or $r=p\pm1$.

In this paper we will extend the result of Kimball and Webb to
arbitrary Lucas sequences. For any positive integer $n$, let $w_n$
be the largest divisor of $u_n$ prime to $u_1,u_2,\ldots,u_{n-1}$.
Here $w_n$ was firstly introduced by Hu and Sun [2] in an
extension of the Lucas congruence for Lucas' $u$-nomial
coefficients.
\proclaim{Theorem 1.1} Let $n\geq 5$ be a positive
integer. Then
$$
\sum_{j=1}^{n-1}\frac{v_j}{u_j}\equiv\frac{(n^2-1)\Delta}{6}\cdot\frac{u_n}{v_n}\pmod{w_n^2}.\tag 1.3
$$
\endproclaim
It is easy to check that either all $u_n$ are odd when $n\geq 1$, or one of $u_2=A$ and $u_3=A^2-B$ is even.
So $w_n$ is odd for any $n>3$. Also we can verify that either $u_n$ is prime to $3$ for each $n\geq 1$, or
$3$ divides one of $u_2$, $u_3$ and $u_4=A^3-2AB$. Hence $3\nmid w_n$ provided that $n>4$.
Finally we mention that $w_n$ is always prime to $v_n$ when $n\geq 3$. Indeed, since
$$
u_n=Au_{n-1}-Bu_{n-2}\text{ and }(w_n, Au_{n-1})=(w_n, u_2u_{n-1})=1,
$$
we have $w_n$ and $B$ are co-prime. And from
$$
v_n=u_{n+1}-Bu_{n-1}=Au_n-2Bu_{n-1},
$$
it follows that $(w_n, v_n)=(w_n, 2Bu_{n-1})=1$.

The Fibonacci numbers $F_0, F_1,\ldots$ are given by
$$
F_0=0,\ F_1=1\text{ and }F_n=F_{n-1}+F_{n-2}\text{ for }n\geq 2.
$$
And the Lucas numbers $L_0, L_1,\ldots$ are given by
$$
L_0=2,\ L_1=1\text{ and }L_n=L_{n-1}+L_{n-2}\text{ for }n\geq 2.
$$
Then by Theorem 1.1, we immediately have
\proclaim{Corollary 1.2} Let $p\geq 5$ be a prime. Let $n$ be the least positive
integer such that $p\mid F_n$. Then
we have
$$
\align
\sum_{j=1}^{n-1}\frac{L_j}{F_j}\equiv&\frac{5(n^2-1)}{6}\cdot\frac{F_n}{L_n}\pmod{p^2}.\tag 1.4
\endalign
$$
\endproclaim
The proof of Theorem 1.1 will be given in the next section.
\heading
2. Proof of Theorem 1.1
\endheading
For any $n\in\N$, the $q$-integer $[n]_q$ is given by
$$
[n]_q=\frac{1-q^n}{1-q}=1+q+\cdots+q^{n-1}.
$$
Now we consider $[n]_q$ as the polynomial in the variable $q$. Recenty Shi and Pan [5] established a $q$-analogue of (1.1)
for prime $p\geq 5$:
$$
\sum_{j=1}^{p-1}\frac{1}{[j]_q}\equiv
\frac{p-1}{2}(1-q)+\frac{p^2-1}{24}(1-q)^2[p]_q\pmod{[p]_q^2}.\tag
2.1
$$

Let
$$
\Phi_n(q)=\prod_{\Sb 1\leq k<n\\(k,n)=1\endSb}(q-\zeta_n^{k})
$$
be the $n$-th cyclotomic polynomial, where $\zeta_n=e^{2\pi i/n}$.
We know that $\Phi_n(q)$ is a polynomial with integral
coefficients, and clearly $\Phi_n(q)$ is prime to $[j]_q$ for each
$1\leq j<n$. Indeed, using similar method, we can easily extend
(2.1) as follows: \proclaim{Proposition 2.1} Let $n$ be a positive
integer. Then
$$
24\sum_{j=1}^{n-1}\frac{1}{[j]_q}\equiv12(n-1)(1-q)+(n^2-1)(1-q)^2[n]_q\pmod{\Phi_n(q)^2}.\tag 2.2
$$
\endproclaim
For the proof of (2.1) and (2.2), the reader may refer to [5].
From (2.2), we deduce that
$$
\align
12\sum_{j=1}^{n-1}\frac{1+q^j}{[j]_q}=&12\sum_{j=1}^{n-1}\frac{2-(1-q^j)}{[j]_q}\\
=&24\sum_{j=1}^{n-1}\frac{1}{[j]_q}-12(n-1)(1-q)\\
\equiv&(n^2-1)(1-q)^2[n]_q\pmod{\Phi_n(q)^2}.
\endalign
$$
And the above congruence can be rewritten as
$$
\bigg(12\sum_{j=1}^{n-1}\frac{1+q^j}{[j]_q}-(n^2-1)(1-q)(1-q^n)\bigg)\prod_{j=1}^{n-1}[j]_q\equiv 0\pmod{\Phi_n(q)^2}.
$$
Since $\Phi_n(q)$ is a primitive polynomial, by Gauss's lemma (cf.
[4], Chapter IV Theorem 2.1 and Corollary 2.2), there exists a
polynomial $G(q)$ with integral coefficients such that
$$
\bigg(12\sum_{j=1}^{n-1}\frac{1+q^j}{[j]_q}-(n^2-1)(1-q)(1-q^n)\bigg)\prod_{j=1}^{n-1}[j]_q=G(q)\Phi_n(q)^2.\tag 2.3
$$

{\noindent{\it Proof of Theorem 1.1}}. When $\Delta=0$, the
theorem reduces to Wolstenholme's congruence (1.1). So below we
assume that $\Delta\not=0$, i.e., $\alpha\not=\beta$. Let $p$ be a
prime with $p\mid w_n$, and let $m$ be the integer such that
$p^m\mid w_n$ but $p^{m+1}\nmid w_n$. Obviously we only need to
show that
$$
\sum_{j=1}^{n-1}\frac{v_j}{u_j}\equiv\frac{(n^2-1)\Delta}{6}\cdot\frac{u_n}{v_n}\pmod{p^{2m}}
$$
for each such $p$ and $m$.

Let $\K=\Q(\sqrt{\Delta})$ and let $\CO_{\K}$ be the ring of algebraic integers in $\K$. Clearly $\alpha, \beta\in\CO_{\K}$.
Let $(p)$ denote the ideal generated by $p$ in $\CO_{\K}$.
We know that if $\jacob{\Delta}{p}=-1$ then
$(p)$ is prime in $\CO_{\K}$, where $\jacob{\cdot}{p}$ is the Legendre symbol. Also there exist two distinct prime ideal $\fp$
and $\fp'$ such that $(p)=\fp\fp'$ provided that $\jacob{\Delta}{p}=1$. Finally when $p\mid\Delta$, $(p)$ is the square of a prime ideal $\fp$. The reader can find the details in [1].
Let
$$
\fP=\cases(p)\qquad&\text{ if }\jacob{\Delta}{p}=-1\text{ or }0,\\\fp\qquad&\text{ if }\jacob{\Delta}{p}=1.\endcases
$$
Obviously either $\alpha$ or $\beta$ is prime to $\fP$, otherwise
we must have $\fP$ is not prime to $u_j$ for any $j\geq 2$, which
implies that $p\mid u_j$. Without loss of generality, we may
assume that $\beta$ is prime to $\fP$.
\proclaim{Lemma 2.2} Let
$p$ be a prime and $k\in\Z$. Suppose that $\jacob{\Delta}{p}=1$
and $(p)=\fp\fp'$. Then for any $m\in\Z^+$, $\fp^m\mid k$ implies
that $p^m\mid k$.
\endproclaim
\Proof. Observe that $\sigma:\,\sqrt{\Delta}\longmapsto-\sqrt{\Delta}$ is an automorphism over $\K$. Also we know that $\sigma(\fp)=\fp'$. Hence
$$
{\fp'}^m=\sigma(\fp^m)\mid\sigma(k)=k.
$$
Since $\fp$ and $\fp'$ are distinct prime ideals, by the unique factorization theorem, we have
$(p)^m=\fp^m{\fp'}^m$ divides $k$.\qed

Now it suffices to prove that
$$
\sum_{j=1}^{n-1}\frac{v_j}{u_j}\equiv\frac{(n^2-1)\Delta}{6}\cdot\frac{u_n}{v_n}\pmod{\fP^{2m}}.
$$
For any $l\in\Z^{+}$, let
$$
\Phi_l(\alpha,\beta)=\beta^{\varphi(l)}\Phi_l(\alpha/\beta)=\prod_{\Sb
1\leq d\leq l\\(d,l)=1\endSb}(\alpha-\zeta_l^{d}\beta),
$$
where $\varphi$ is the Euler totient function. Apparently
$\Phi_l(\alpha,\beta)\in\CO_{\K}$. Notice that
$$
u_l=\frac{\alpha^l-\beta^l}{\alpha-\beta}=\prod_{1<d, d\mid l}\beta^{\varphi(d)}\Phi_d(\alpha/\beta)=\prod_{1<d, d\mid l}\Phi_d(\alpha,\beta).
$$
Hence $u_l$ is always divisible by $\Phi_l(\alpha,\beta)$.
Then we have $w_n$ divides
$$
\Phi_n(\alpha,\beta)=\frac{u_n}{\prod_{\Sb 1<d<n\\d\mid n\endSb}\Phi_d(\alpha,\beta)}
$$
since $w_n$ is prime to $u_d$ whenever $1\leq d<n$.

Substituting $\alpha/\beta$ for $q$ in (2.3), and noting that
$$
u_j=\beta^{j-1}\frac{1-(\alpha/\beta)^j}{1-\alpha/\beta}=\beta^{j-1}[j]_{\alpha/\beta}\text{ and }v_j=\beta^{j}(1+(\alpha/\beta)^j),
$$
we obtain that
$$
\bigg(12\beta^{-1}\sum_{j=1}^{n-1}\frac{v_j}{u_j}-(n^2-1)\beta^{-n-1}(\alpha-\beta)^2u_n\bigg)\prod_{j=1}^{n-1}\beta^{1-j}u_j=\beta^{-2\varphi(n)}G(\alpha/\beta)\Phi_d(\alpha,\beta)^2.
$$
As $(w_n,6)=1$ and $\fP$ is prime to $\beta$,
we conclude that
$$
\sum_{j=1}^{n-1}\frac{v_j}{u_j}-\frac{(n^2-1)\Delta}{12}\beta^{-n}u_n\equiv0\pmod{(\fP^m)^2}.
$$
Finally, since
$$
\alpha^n=\frac{1}{2}(v_n+u_n\sqrt{\Delta})\text{ and }\beta^n=\frac{1}{2}(v_n-u_n\sqrt{\Delta}),
$$
we have
$$
\alpha^n\equiv\beta^n\equiv v_n/2\pmod{w_n}.
$$
All are done.\qed

\Ack. I thank the referee for the valuable comments on this paper. I also thank my advisor, Professor Zhi-Wei Sun, for his helpful suggestions.

\enddocument